\documentclass[12pt]{amsart}
\usepackage{amssymb}
\usepackage{amsfonts}
\usepackage{color,soul}

\newcommand{\nc}{\newcommand}
\nc{\thref}[1]{Theorem~\ref{theo:#1}}
\nc{\selabel}[1]{\label{sect:#1}}
\nc{\seref}[1]{Section~\ref{sect:#1}}
\nc{\lelabel}[1]{\label{lemm:#1}}
\nc{\leref}[1]{Lemma~\ref{lemm:#1}}
\nc{\prlabel}[1]{\label{prop:#1}}
\nc{\prref}[1]{Proposition~\ref{prop:#1}}
\nc{\colabel}[1]{\label{coro:#1}}
\nc{\coref}[1]{Corollary~\ref{coro:#1}}
\nc{\exlabel}[1]{\label{exam:#1}}
\nc{\exref}[1]{Example~\ref{exam:#1}}
\nc{\delabel}[1]{\label{defi:#1}}
\nc{\deref}[1]{Definition~\ref{defi:#1}}
\nc{\eqlabel}[1]{\label{equa:#1}}
\nc{\relabel}[1]{\label{rema:#1}}
\nc{\reref}[1]{Lemma~\ref{rema:#1}}
\providecommand{\operatorname}[1]{\mathrm{#1}\,}
\nc{\Hom}{\operatorname{Hom}} \nc{\Mor}{\operatorname{Mor}}
\nc{\Aut}{\operatorname{Aut}} \nc{\Ann}{\operatorname{Ann}}
\nc{\Ker}{\operatorname{Ker}} \nc{\Trace}{\operatorname{Trace}}
\nc{\Char}{\operatorname{Char}} \nc{\Mod}{\operatorname{Mod}}
\nc{\End}{\operatorname{End}} \nc{\Spec}{\operatorname{Spec}}
\nc{\Span}{\operatorname{Span}} \nc{\sgn}{\operatorname{sgn}}
\nc{\Id}{\operatorname{Id}} \nc{\Com}{\operatorname{Com}}
\nc{\rank}{\operatorname{rank}}
\nc{\Clausen}{\operatorname{Cl}}
\nc{\Li}{\operatorname{Li}}
\nc{\Ls}{\operatorname{Ls}}

\let\:=\colon

\newtheorem{de}{Definition}[section]
\newtheorem{lm}[de]{Lemma}
\newtheorem{pr}[de]{Proposition}
\newtheorem{co}[de]{Corollary}
\newtheorem{re}[de]{Remark}
\newtheorem{res}[de]{Remarks}
\newtheorem{te}[de]{Theorem}
\newtheorem{ex}[de]{Example}
\newtheorem{exs}[de]{Examples}

\def\bex{\begin{ex}}
\def\eex{\end{ex}}
\def\bexs{\begin{exs}}
\def\eexs{\end{exs}}
\def\bl{\begin{lm}}
\def\el{\end{lm}}
\def\bc{\begin{co}}
\def\ec{\end{co}}
\def\bt{\begin{te}}
\def\et{\end{te}}
\def\bpr{\begin{pr}}
\def\epr{\end{pr}}
\def\br{\begin{re}}
\def\er{\end{re}}
\def\brs{\begin{res}}
\def\ers{\end{res}}
\def\bd{\begin{de}}
\def\ed{\end{de}}
\def\be{\begin{equation}}
\def\ee{\end{equation}}
\def\bea{\begin{eqnarray*}}
\def\eea{\end{eqnarray*}}
\def\bp{\begin{proof}}
\def\ep{\end{proof}}

\def\qed{\hfill\Box}
\def\NN{{\mathbb N}}

\let\:=\colon

\begin{document}

\title[Generalized Log-sine integrals and Bell polynomials]{Generalized Log-sine integrals and Bell polynomials}

\begin{abstract}
In this paper, we investigate the integral of $x^n\log^m(\sin(x))$ for natural numbers $m$ and $n$. In doing so, we recover some well-known results and remark on some relations to the log-sine integral $\Ls_{n+m+1}^{(n)}(\theta)$. Later, we use properties of Bell polynomials to find a closed expression for the derivative of the central binomial and shifted central binomial coefficients in terms of polygamma functions and harmonic numbers.
\end{abstract}

\author{Derek Orr}

\thanks{2010 \textit{Mathematics Subject Classification}. Primary 33E20, 11B73. Secondary 11M32.}

\keywords{Log-sine integral, Riemann zeta function, Bell polynomial, harmonic numbers, Euler sum, binomial coefficients}

\maketitle

\section{Introduction and Preliminaries}

\vspace{0.3cm}

The functions

$$ \Ls_{n}(\theta) : = -\int_{0}^{\theta} \Big(\log\Big|2\sin\Big(\frac{x}{2}\Big)\Big|\Big)^{n-1} \hspace{3pt} dx $$

\vspace{0.3cm}

and 

$$ \Ls_{n}^{(m)}(\theta) := -\int_{0}^{\theta} x^m \Big(\log\Big|2\sin\Big(\frac{x}{2}\Big)\Big|\Big)^{n-m-1} \hspace{3pt} dx $$

\vspace{0.3cm}

have been widely studied in previous papers (see \cite{Borwein}, \cite{Borwein1}, \cite{Bowman}, \cite{Choi}, \cite{Choi0}, \cite{DaK}, \cite{NanYue}). A very nice identity was given in \cite{Borwein} by expressing

$$ S(k) := \sum_{n=1}^{\infty} \frac{1}{n^k\binom{2n}{n}} $$

as

$$ S(k) = \frac{(-2)^{k-1}}{(k-2)!}\Ls_{k}^{(1)}\Big(\frac{\pi}{3}\Big), k \in \NN.$$

\vspace{0.5cm}

Here, we will focus on a similar integral,

$$F(n,m,z) = \int_{0}^{z} x^n\sin^{2m}(x) \hspace{3pt} dx. $$

\vspace{0.3cm}

Further, we can define 

$$G(n,m,z) = \int_{0}^{z} x^n\Big(2\sin\Big(\frac{x}{2}\Big)\Big)^{2m} \hspace{3pt} dx, $$

\vspace{0.3cm}

and we can easily see that

\begin{equation}
G(n,m,2z) = 4^m2^{n+1}F(n,m,z),
\end{equation}

\vspace{0.3cm}

and

\begin{equation}
-\frac{\partial^p G(n,m,z)}{\partial m^p}\bigg|_{(n,0,z)} = 2^p\Ls_{p+n+1}^{(n)}(z).
\end{equation}

\vspace{0.4cm}

As we discuss the behavior of $F(n,m,z)$, we will add in remarks for $G(n,m,z)$ and thus for $\Ls_{p+n+1}^{(n)}(z)$. Next, we introduce the Riemann zeta function and the polylogarithm function.

\begin{equation}
\displaystyle\zeta(s):=  \left\{
\begin{array}{ll}
      \displaystyle\sum_{n=1}^{\infty}\frac{1}{n^s}=\frac{1}{1-2^{-s}}\sum_{n=1}^{\infty}\frac{1}{(2n-1)^s}, & \operatorname{Re}(s)>1, \\ \\
      \displaystyle\frac{1}{1-2^{1-s}}\sum_{n=1}^{\infty}\frac{(-1)^{n-1}}{n^s},  & \operatorname{Re}(s)>0, \hspace{3pt} s\neq 1,\\
\end{array} 
\right.
\end{equation}

\vspace{0.3cm}

and

\begin{equation}
\Li_{n}(z) = \sum_{k=1}^{\infty} \frac{z^k}{k^n}, \hspace{5pt} n \in \NN\backslash \{1\}, \hspace{5pt} |z| \leq 1.
\end{equation}

\vspace{0.3cm}

Euler discovered the now famous closed formula for $\zeta(2n)$, given by

\begin{equation} 
\zeta(2k)= \sum_{n=1}^{\infty} \frac{1}{n^{2k}} = \frac{(-1)^{k+1}B_{2k}(2\pi)^{2k}}{2(2k)!}, \hspace{5pt} k \in \mathbb{N}_{0},
\end{equation}

\vspace{0.3cm}

where $B_{n}$ are the Bernoulli numbers, defined by

$$\frac{z}{e^z-1}=\sum_{n=0}^{\infty}\frac{B_{n}}{n!}z^n, \hspace{5pt} |z|<2\pi.$$

\vspace{0.2cm}

It is clear that

\begin{equation}
\Li_{n}(1) = \zeta(n), \hspace{10pt} \Li_{n}(-1) = -(1-2^{1-n})\zeta(n), \hspace{5pt} n \in \NN \backslash \{1\}.
\end{equation}

\vspace{0.4cm}

We also introduce the generalized hypergeometric function

$$ {}_{p}F_{q}(a_{1},a_{2},\dots,a_{p};b_{1},b_{2},\dots,b_{q};z) = \sum_{k=0}^{\infty} \frac{(a_{1})_{k}(a_{2})_{k}\dots(a_{p})_{k}}{(b_{1})_{k}(b_{2})_{k}\dots(b_{q})_{k}}\frac{z^k}{k!},$$

\vspace{0.3cm}

where

\begin{equation}
(a)_{k} = \frac{\Gamma(a+k)}{\Gamma(a)} = a(a+1)(a+2)\dots(a+k-1)
\end{equation}

\vspace{0.3cm}

is the Pochhammer symbol or rising factorial. If $a_{1} = a_{2} = \dots = a_{i} = a$, we will use the notation ${}_{p}F_{q}(\{a\}^i,a_{i+1},\dots,a_{p};b_{1},b_{2},\dots,b_{q};z)$. A special case used in the paper is

$${}_{q+1}F_{q}(\{1\}^{q+1};\{2\}^{q};z) = \sum_{k=0}^{\infty} \frac{(1)_{k}(1)_{k}\dots(1)_{k}}{(2)_{k}(2)_{k}\dots(2)_{k}}\frac{(1)_{k}z^k}{k!}, $$

\vspace{0.3cm}

which becomes

\begin{equation}
{}_{q+1}F_{q}(\{1\}^{q+1};\{2\}^{q};z) = \sum_{k=0}^{\infty} \frac{z^k}{(k+1)^q} = \frac{1}{z}\Li_{q}(z).
\end{equation}

\vspace{0.4cm}

Since this paper will involve the derivative of the gamma function, we define polygamma function

$$ \psi^{(n)}(z) := \frac{d^{n+1}}{dz^{n+1}} \log\Gamma(z), \hspace{5pt} n \in \NN_{0}.$$

\vspace{0.3cm}

The reflection and recursive formulas are given by

\begin{equation}
\psi^{(n)}(1-z) = (-1)^n\Big(\psi^{(n)}(z) + \pi\frac{d^n}{dz^n}(\cot(\pi z))\Big),
\end{equation}

\vspace{0.3cm}

and

\begin{equation}
\psi^{(n)}(1+z) = \psi^{(n)}(z) + \frac{(-1)^nn!}{z^{n+1}},
\end{equation}

\vspace{0.3cm}

respectively. When $z=1$, we have

\begin{equation}
\psi^{(n)}(1) = (-1)^{n+1}n!\zeta(n+1),
\end{equation}

\vspace{0.3cm}

and in general,

$$ \psi^{(n)}(z) = (-1)^{n+1}n!\big(\zeta(n+1)-H_{z-1}^{(n+1)}\big), $$

\vspace{0.3cm}

where $H_{m}^{(n)}$ are the generalized harmonic numbers. When $m=0$, it is understood that $H_{0}^{(n)}=0$. Lastly, we introduce the multiple zeta function

$$\zeta(s_{1},s_{2},\dots,s_{k}) = \sum_{n_{1}>n_{2}>\dots>n_{k}>0} \frac{1}{n_{1}^{s_{1}}\dots n_{k}^{s_{k}}}. $$

\vspace{0.3cm}

If $s_{1}=s_{2}=\dots=s_{i}=s$, it is common to denote the multiple zeta function as $\zeta(\{s\}^i,s_{i+1},s_{i+2},\dots,s_{k})$. Further, a horizontal bar will be given to a variable if its sum is alternating. For example,

$$\zeta(2,1,\overline{4}) = \sum_{n_{1}>n_{2}>n_{3}>0} \frac{(-1)^{n_3}}{n_{1}^2n_{2}n_{3}^4}.$$

\vspace{0.3cm}

In particular,

$$\zeta(\overline{2j+1},1) = \sum_{n_{1}>n_{2}>0} \frac{(-1)^{n_{1}}}{n_{1}^{2j+1}n_{2}} = \sum_{n_{1}=2}^{\infty} \frac{(-1)^{n_{1}}}{n_{1}^{2j+1}}\sum_{n_{2}=1}^{n_{1}-1}\frac{1}{n_{2}} = \sum_{n_{1}=2}^{\infty} \frac{(-1)^{n_{1}}H_{n_{1}-1}}{n_{1}^{2j+1}}, $$

\vspace{0.4cm}

and using $H_{k-1} = H_{k}-1/k$ and rearranging, we have the formula

$$ \sum_{k=1}^{\infty} \frac{(-1)^kH_{k}}{k^{2j+1}} = \zeta(\overline{2j+1},1)+\sum_{k=1}^{\infty} \frac{(-1)^k}{k^{2j+2}} = \zeta(\overline{2j+1},1)-(1-2^{-2j-1})\zeta(2j+2). $$

\vspace{0.3cm}

Another famous formula for harmonic sums studied in \cite{Adamchik} is

$$ \sum_{k=1}^{\infty} \frac{H_{k}}{k^n} = \frac{1}{2}(n+2)\zeta(n+1)-\frac{1}{2}\sum_{k=1}^{n-2}\zeta(k+1)\zeta(n-k). $$

\vspace{0.3cm}

The multiple zeta function, as well as the other functions mentioned, have been studied and each has a wide variety of applications in mathematics and physics (see \cite{DaK}, \cite{DaK1}, \cite{Koyama}, \cite{Sofo}, \cite{Yost}). In this paper, we will find a formula for the partial derivatives of $F(n,m,z)$ with respect to $m$ for specific $z$ and hence a formula for $\Ls_{n}^{(m)}(\theta)$ in terms of derivatives of binomial coefficients. In the latter half of the paper, using Bell polynomials, we give explicit formulas for these derivatives in terms of harmonic numbers and polygamma functions. We begin by introducing some equations that can be found in \cite{GradRyz}. We have

\begin{multline}
F(n,m,z) = \int_{0}^{z} x^n\sin^{2m}(x) \hspace{3pt} dx = \frac{\binom{2m}{m}z^{n+1}}{4^m(n+1)}\\+\sum_{k=0}^{m-1}\frac{(-1)^{m+k}}{2^{2m-1}}\binom{2m}{k}\int_{0}^{z} x^n\cos(2(m-k)x) \hspace{3pt} dx,
\end{multline}

\vspace{0.3cm}

and

\begin{equation}
\int_{0}^{z} x^n\cos(ax) \hspace{3pt} dx = \sum_{j=0}^{n}\frac{n!z^{n-j}}{(n-j)!a^{j+1}}\sin\Big(az+\frac{\pi}{2}j\Big)-\frac{n!}{a^{n+1}}\sin\Big(\frac{\pi n}{2}\Big).
\end{equation}

\vspace{0.3cm}

Combining them and reindexing the sum on $k$, we have

\begin{multline}
F(n,m,z) = \int_{0}^{z} x^n\sin^{2m}(x) \hspace{3pt} dx = \frac{1}{4^m}\binom{2m}{m}\frac{z^{n+1}}{n+1}\\+n!\sum_{k=1}^{\infty}\frac{(-1)^k}{2^{2m-1}}\binom{2m}{m+k}\bigg(\sum_{j=0}^{n}\frac{z^{n-j}\sin(2kz+\pi j/2)}{(2k)^{j+1}(n-j)!}-\frac{\sin(\pi n/2)}{(2k)^{n+1}}\bigg).
\end{multline}

\vspace{0.4cm}

Using $(1)$, we see

\begin{multline}
G(n,m,2z) = \int_{0}^{2z} x^n\Big(2\sin\Big(\frac{x}{2}\Big)\Big)^{2m} \hspace{3pt} dx = \binom{2m}{m}\frac{(2z)^{n+1}}{n+1}\\+2n!\sum_{k=1}^{\infty}(-1)^k\binom{2m}{m+k}\bigg(\sum_{j=0}^{n}\frac{(2z)^{n-j}\sin(2kz+\pi j/2)}{k^{j+1}(n-j)!}-\frac{\sin(\pi n/2)}{k^{n+1}}\bigg).
\end{multline}

\vspace{0.4cm}

Our last introductory remark brings us back to the generalized hypergeometric function. First, using $(7)$, $\displaystyle \frac{(m-k)_{k}}{(m+2)_{k}} = \frac{\binom{2m}{m+k+1}}{\binom{2m}{m+1}}$. With this, we see

$${}_{s+2}F_{s+1}(\{1\}^{s+1},1-m;\{2\}^{s},m+2;-z) = \sum_{k=0}^{\infty} \frac{(1)_{k}(1)_{k}\dots(1)_{k}}{(2)_{k}(2)_{k}\dots(2)_{k}}\frac{(1-m)_{k}}{(m+2)_{k}}\frac{(1)_{k}(-z)^k}{k!}$$

$$= \sum_{k=0}^{\infty} \frac{(-1)^k(m-k)_{k}(-z)^k}{(k+1)^s(m+2)_{k}} = \sum_{k=0}^{\infty} \frac{z^k\binom{2m}{m+k+1}}{(k+1)^s\binom{2m}{m+1}},$$

\vspace{0.3cm}

and lastly, changing the index of our sum, we have

\begin{equation}
\sum_{k=1}^{\infty} \frac{z^{k-1}}{k^s}\binom{2m}{m+k} = \binom{2m}{m+1} {}_{s+2}F_{s+1}(\{1\}^{s+1},1-m;\{2\}^{s},m+2;-z).
\end{equation}

\vspace{0.4cm}

\section{Log-sine integral for $z=2\pi$}

\bt For $n, m \in \mathbb{N}_{0}$,

\begin{equation}
\int_{0}^{\pi} x^n\sin^{2m}(x) \hspace{3pt} dx  = \frac{\pi^n}{4^m}\bigg(\frac{\pi\binom{2m}{m}}{n+1}-n!\sum_{j=1}^{\lfloor \frac{n}{2} \rfloor}\frac{(-1)^j}{(2\pi)^{2j-1}(n+1-2j)!}\sum_{k=1}^{\infty}\frac{(-1)^k}{k^{2j}}\binom{2m}{m+k}\bigg).
\end{equation}

\et

\vspace{0.4cm}

\textit{Proof.} Letting $z=\pi$ in $(12)$,

\begin{multline*}
F(n,m,\pi) = \int_{0}^{\pi} x^n\sin^{2m}(x) \hspace{3pt} dx = \frac{1}{4^m}\binom{2m}{m}\frac{\pi^{n+1}}{n+1}\\+n!\sum_{k=1}^{\infty}\frac{(-1)^k}{2^{2m-1}}\binom{2m}{m+k}\bigg(\sum_{j=0}^{n}\frac{\pi^{n-j}\sin(2k\pi+\pi j/2)}{(2k)^{j+1}(n-j)!}-\frac{\sin(\pi n/2)}{(2k)^{n+1}}\bigg)
\end{multline*}

$$= \frac{1}{4^m}\binom{2m}{m}\frac{\pi^{n+1}}{n+1}+n!\sum_{k=1}^{\infty}\frac{(-1)^k}{2^{2m-1}}\binom{2m}{m+k}\sum_{j=0}^{n-1}\frac{\pi^{n-j}\sin(2k\pi+\pi j/2)}{(2k)^{j+1}(n-j)!} $$

$$= \frac{\pi^n}{4^m}\bigg(\binom{2m}{m}\frac{\pi}{n+1}+n!\sum_{j=0}^{n-1}\sum_{k=1}^{\infty}\frac{(-1)^k}{k^{j+1}}\binom{2m}{m+k}\frac{\sin(\pi j/2)}{(2\pi)^j(n-j)!}\bigg). $$

\vspace{0.4cm}

Using $\sin(\pi j/2) = \delta_{\lfloor \frac{j-1}{2} \rfloor, \frac{j-1}{2}}(-1)^{\frac{j-1}{2}}$ we have,

$$
\int_{0}^{\pi} x^n\sin^{2m}(x) \hspace{3pt} dx  = \frac{\pi^n}{4^m}\bigg(\frac{\pi\binom{2m}{m}}{n+1}+n!\sum_{j=0}^{\lfloor \frac{n-2}{2} \rfloor}\frac{(-1)^j}{(2\pi)^{2j+1}(n-2j-1)!}\sum_{k=1}^{\infty}\frac{(-1)^k}{k^{2j+2}}\binom{2m}{m+k}\bigg),$$

\vspace{0.3cm}

and changing the index on $j$, the proof is complete. $\qed$ 

\vspace{0.5cm}

Note from this that

\begin{equation}
G(n,m,2\pi) = 2(2\pi)^n\bigg(\frac{\pi\binom{2m}{m}}{n+1}-n!\sum_{j=1}^{\lfloor \frac{n}{2} \rfloor}\frac{(-1)^j}{(2\pi)^{2j-1}(n+1-2j)!}\sum_{k=1}^{\infty}\frac{(-1)^k}{k^{2j}}\binom{2m}{m+k}\bigg).
\end{equation}

\vspace{0.4cm}

Now taking the derivative of $(17)$ and using $(16)$, we find

\begin{multline*}
\frac{\partial F}{\partial m}\bigg|_{(n,0,\pi)} = 2\int_{0}^{\pi} x^n\log(\sin(x)) \hspace{3pt} dx = \pi^n\bigg( \frac{\pi}{n+1}\frac{d}{dm}\bigg(\frac{\binom{2m}{m}}{4^m}\bigg)\bigg|_{m=0}\\+ \frac{d}{dm}\binom{2m}{m+1}\bigg|_{m=0} \sum_{j=1}^{\lfloor \frac{n}{2} \rfloor}\frac{n!(-1)^j{}_{2j+2}F_{2j+1}(\{1\}^{2j+2};\{2\}^{2j+1};1)}{(n+1-2j)!(2\pi)^{2j-1}}\bigg)
\end{multline*}

$$ = \pi^n\bigg( \frac{-2\pi\log2}{n+1}+\sum_{j=1}^{\lfloor \frac{n}{2} \rfloor}\frac{n!(-1)^j\zeta(2j+1)}{(n+1-2j)!(2\pi)^{2j-1}}\bigg) = -2\pi^{n+1}\bigg( \frac{\log2}{n+1}-\sum_{j=1}^{\lfloor \frac{n}{2} \rfloor}\frac{n!(-1)^j\zeta(2j+1)}{(n+1-2j)!(2\pi)^{2j}}\bigg), $$

\vspace{0.4cm}

where $(6)$ and $(8)$ have been used. Taking more partial derivatives, we see that

\begin{multline}
\frac{\partial^p F}{\partial m^p}\bigg|_{(n,0,\pi)} = 2^p\int_{0}^{\pi} x^n\log^p(\sin(x)) \hspace{3pt} dx = \pi^n\bigg(\frac{\pi}{n+1}\frac{d^p}{dm^p}\bigg(\frac{\binom{2m}{m}}{4^m}\bigg)\bigg|_{m=0}\\-\sum_{k=1}^{\infty} \frac{\partial^p}{\partial m^p}\bigg(\frac{\binom{2m}{m+k}}{4^m}\bigg)\bigg|_{m=0}\sum_{j=1}^{\lfloor \frac{n}{2} \rfloor}\frac{n!(-1)^{j+k}}{(n+1-2j)!(2\pi)^{2j-1}k^{2j}}\bigg). 
\end{multline}

%\vspace{0.3cm}

%Clearly, if $n \in \{0,1\}$,

%\begin{equation}
%2^p\int_{0}^{\pi} x^n\log^p(\sin(x)) \hspace{3pt} dx = \frac{\pi^{n+1}}{n+1}\frac{d^p}{dm^p}\bigg(\frac{\binom{2m}{m}}{4^m}\bigg)\bigg|_{m=0}. 
%\end{equation}

\vspace{0.4cm}

Below we compute a few integrals for $p>1$ and $n>0$.

$$ \int_{0}^{\pi} x\log^2(\sin(x)) \hspace{3pt} dx = \frac{\pi^4}{24}+\frac{\pi^2}{2}\log^22$$

$$ \int_{0}^{\pi} x^2\log^2(\sin(x)) \hspace{3pt} dx = \frac{13\pi^5}{360}+\pi\zeta(3)\log2+\frac{\pi^3}{3}\log^22$$

$$ \int_{0}^{\pi} x^3\log^2(\sin(x)) \hspace{3pt} dx = \frac{\pi^6}{30}+\frac{3\pi^2}{2}\zeta(3)\log2+\frac{\pi^4}{4}\log^22$$

$$ \int_{0}^{\pi} x^4\log^2(\sin(x)) \hspace{3pt} dx = \frac{37\pi^7}{1260}+2\pi^3\zeta(3)\log2-3\pi\zeta(5)\log2+\frac{3\pi}{2}\zeta^2(3)+\frac{\pi^5}{5}\log^22$$

%$$ \int_{0}^{\pi} x^5\log^2(\sin(x)) \hspace{3pt} dx = \frac{25\pi^8}{1008}+\frac{5\pi^4}{2}\zeta(3)\log2-\frac{15\pi^2}{2}\zeta(5)\log2+\frac{15\pi^2}{4}\zeta^2(3)+\frac{\pi^6}{6}\log^22 $$

$$ \int_{0}^{\pi} x\log^3(\sin(x)) \hspace{3pt} dx =\frac{3\pi^2}{8}\zeta(3)-\frac{\pi^4}{16}\log2-\frac{\pi^2}{4}\log^32 $$

\vspace{0.5cm}

If we use $(18)$ instead,

\begin{multline}
\frac{\partial^p G}{\partial m^p}\bigg|_{(n,0,2\pi)} = -2^p\Ls_{p+n+1}^{(n)}(2\pi) =  2^p\int_{0}^{2\pi} x^n\log^p\Big(2\sin\Big(\frac{x}{2}\Big)\Big) \hspace{3pt} dx \\= 2(2\pi)^n\bigg(\frac{\pi}{n+1}\frac{d^p}{dm^p}\binom{2m}{m}\bigg|_{m=0}-\sum_{k=1}^{\infty} \frac{\partial^p}{\partial m^p}\binom{2m}{m+k}\bigg|_{m=0}\sum_{j=1}^{\lfloor \frac{n}{2} \rfloor}\frac{n!(-1)^{j+k}}{(n+1-2j)!(2\pi)^{2j-1}k^{2j}}\bigg).
\end{multline}

%\vspace{0.3cm}

%and if $n \in \{0,1\}$,

%\begin{equation}
%-2^p\Ls^{(n)}_{p+n+1}(2\pi) =  2^p\int_{0}^{2\pi} x^n\log^p\Big(2\sin\Big(\frac{x}{2}\Big)\Big) \hspace{3pt} dx = \frac{(2\pi)^{n+1}}{n+1}\frac{d^p}{dm^p}\binom{2m}{m}\bigg|_{m=0}.
%\end{equation}

\vspace{0.5cm}

Below we compute a few of these integrals.

$$ \Ls_{4}^{(1)}(2\pi) = -\frac{\pi^4}{6}, \hspace{10pt} \Ls_{5}^{(1)}(2\pi) = 3\pi^2\zeta(3), \hspace{10pt} \Ls_{5}^{(2)}(2\pi) = -\frac{13\pi^5}{45}, \hspace{10pt} \Ls_{6}^{(3)}(2\pi) = -\frac{8\pi^6}{15} $$

$$ \Ls_{7}^{(4)}(2\pi) = -\frac{296\pi^7}{315}-48\pi\zeta^2(3), \hspace{10pt} \Ls_{8}^{(5)}(2\pi) = -\frac{100\pi^8}{63}-240\pi^2\zeta^2(3)$$

\vspace{0.3cm}

\section{Log-sine integral for $z=\pi$}

\vspace{0.4cm}

\bt{For $n, m \in \mathbb{N}_{0}$,}

\begin{multline}
\int_{0}^{\pi/2} x^n\sin^{2m}(x) \hspace{3pt} dx = \frac{1}{4^m}\Big(\frac{\pi}{2}\Big)^n\bigg(\frac{\pi\binom{2m}{m}}{2(n+1)}-n!\sum_{k=1}^{\infty}\binom{2m}{m+k}\\ * \bigg(\sum_{j=1}^{\lfloor \frac{n+1}{2} \rfloor}\frac{(-1)^j}{\pi^{2j-1}k^{2j}(n+1-2j)!}-\delta_{\lfloor \frac{n+1}{2} \rfloor,\frac{n+1}{2}}\frac{(-1)^k(-1)^{\lfloor \frac{n+1}{2} \rfloor}}{\pi^nk^{n+1}}\bigg)\bigg).
\end{multline}

\et

\vspace{0.3cm}

\textit{Proof.} Letting $z=\pi/2$, $(12)$ becomes

\begin{multline*}
F(n,m,\pi/2) = \int_{0}^{\pi/2} x^n\sin^{2m}(x) \hspace{3pt} dx = \frac{1}{4^m}\binom{2m}{m}\frac{(\pi/2)^{n+1}}{n+1}\\+n!\sum_{k=1}^{\infty}\frac{(-1)^k}{2^{2m-1}}\binom{2m}{m+k}\bigg(\sum_{j=0}^{n}\frac{(\pi/2)^{n-j}\sin(k\pi+\pi j/2)}{(2k)^{j+1}(n-j)!}-\frac{\sin(\pi n/2)}{(2k)^{n+1}}\bigg)
\end{multline*}

$$
= \frac{1}{4^m}\Big(\frac{\pi}{2}\Big)^n\bigg(\binom{2m}{m}\frac{\pi}{2(n+1)}+n!\sum_{k=1}^{\infty}(-1)^k\binom{2m}{m+k}\bigg(\sum_{j=0}^{n}\frac{(-1)^k\sin(\pi j/2)}{\pi^jk^{j+1}(n-j)!}-\frac{\sin(\pi n/2)}{\pi^nk^{n+1}}\bigg)\bigg).$$

\vspace{0.4cm}

Using the same analysis as before on the sine function, 

\begin{multline*}
\int_{0}^{\pi/2} x^n\sin^{2m}(x) \hspace{3pt} dx = \frac{1}{4^m}\Big(\frac{\pi}{2}\Big)^n\bigg(\binom{2m}{m}\frac{\pi}{2(n+1)}\\+n!\sum_{k=1}^{\infty}(-1)^k\binom{2m}{m+k}\bigg(\sum_{j=0}^{\lfloor \frac{n-1}{2} \rfloor}\frac{(-1)^k(-1)^j}{\pi^{2j+1}k^{2j+2}(n-1-2j)!}-\delta_{\lfloor \frac{n-1}{2} \rfloor,\frac{n-1}{2}}\frac{(-1)^{\lfloor \frac{n-1}{2} \rfloor}}{\pi^nk^{n+1}}\bigg)\bigg)
\end{multline*}

\begin{multline*}
= \frac{1}{4^m}\Big(\frac{\pi}{2}\Big)^n\bigg(\frac{\pi\binom{2m}{m}}{2(n+1)}-n!\sum_{k=1}^{\infty}\binom{2m}{m+k}\\ * \bigg(\sum_{j=1}^{\lfloor \frac{n+1}{2} \rfloor}\frac{(-1)^j}{\pi^{2j-1}k^{2j}(n+1-2j)!}-\delta_{\lfloor \frac{n+1}{2} \rfloor,\frac{n+1}{2}}\frac{(-1)^k(-1)^{\lfloor \frac{n+1}{2} \rfloor}}{\pi^nk^{n+1}}\bigg)\bigg),
\end{multline*}

\vspace{0.4cm}

which completes the proof. $\qed$

\vspace{0.5cm}

Using $(1)$,

\begin{multline}
G(n,m,\pi) = 2\pi^n\bigg(\frac{\pi\binom{2m}{m}}{2(n+1)}-n!\sum_{k=1}^{\infty}\binom{2m}{m+k}\\ * \bigg(\sum_{j=1}^{\lfloor \frac{n+1}{2} \rfloor}\frac{(-1)^j}{\pi^{2j-1}k^{2j}(n+1-2j)!}-\delta_{\lfloor \frac{n+1}{2} \rfloor,\frac{n+1}{2}}\frac{(-1)^k(-1)^{\lfloor \frac{n+1}{2} \rfloor}}{\pi^nk^{n+1}}\bigg)\bigg).
\end{multline}

\vspace{0.5cm}

Again, taking the derivative of $(21)$ and using $(16)$,

\begin{multline*} 
\frac{\partial F}{\partial m}\Big|_{(n,0,\pi/2)} = 2\int_{0}^{\pi/2} x^n\log(\sin(x)) \hspace{3pt} dx = \Big(\frac{\pi}{2}\Big)^n\Bigg[\frac{\pi}{n+1}\frac{d}{dm}\bigg(\frac{\binom{2m}{m}}{2^{2m+1}}\bigg)\Bigg|_{m=0}\\-n!\frac{d}{dm}\bigg(\binom{2m}{m+1}\bigg)\Bigg|_{m=0} \bigg(\sum_{j=1}^{\lfloor \frac{n+1}{2} \rfloor} \frac{(-1)^j{}_{2j+2}F_{2j+1}(\{1\}^{2j+2};\{2\}^{2j+1};-1)}{(n+1-2j)!\pi^{2j-1}}\\ + \delta_{\lfloor \frac{n+1}{2} \rfloor, \frac{n+1}{2}}\frac{(-1)^{\lfloor \frac{n+1}{2} \rfloor}{}_{n+3}F_{n+2}(\{1\}^{n+3};\{2\}^{n+2};1)}{\pi^n}\bigg)\Bigg]
\end{multline*}

$$ = -2\Big(\frac{\pi}{2}\Big)^{n+1}\Bigg[\frac{\log2}{n+1}+n!\bigg(\sum_{j=1}^{\lfloor \frac{n+1}{2} \rfloor} \frac{(-1)^j(2^{2j}-1)\zeta(2j+1)}{(n-2j+1)!(2\pi)^{2j}} + \delta_{\lfloor \frac{n+1}{2} \rfloor, \frac{n+1}{2}}\frac{(-1)^{\lfloor \frac{n+1}{2} \rfloor}\zeta(n+2)}{\pi^{n+1}}\bigg)\Bigg],$$

\vspace{0.3cm}

where $(6)$ and $(8)$ have been used again. Note this formula is the same as in other papers as well (see \cite{Choi}, \cite{Orr}). Again, taking more partial derivatives, we will have

\begin{multline}
\frac{\partial^p F}{\partial m^p}\bigg|_{(n,0,\pi/2)} = 2^p\int_{0}^{\pi/2} x^n\log^{p}(\sin(x)) \hspace{3pt} dx = \Big(\frac{\pi}{2}\Big)^n\Bigg[ \frac{\pi}{2(n+1)}\frac{d^p}{dm^p}\bigg(\frac{\binom{2m}{m}}{4^m}\bigg)\bigg|_{m=0}\\-n!\sum_{k=1}^{\infty} \frac{\partial^p}{\partial m^p}\bigg(\frac{\binom{2m}{m+k}}{4^m}\bigg)\bigg|_{m=0}\bigg(\sum_{j=1}^{\lfloor \frac{n+1}{2} \rfloor}\frac{(-1)^j}{(n+1-2j)!\pi^{2j-1}k^{2j}}-\delta_{\lfloor \frac{n+1}{2} \rfloor, \frac{n+1}{2}} \frac{(-1)^k(-1)^{\lfloor \frac{n+1}{2}\rfloor}}{\pi^nk^{n+1}}\bigg)\Bigg].
\end{multline}

\vspace{0.4cm}

%Then dividing by $2^p$ will yield the integral value. Surely, if $n=0$, we have

%\begin{equation}
%2^p\int_{0}^{\pi/2} \log^{p}(\sin(x)) \hspace{3pt} dx = \frac{\pi}{2}\frac{d^p}{dm^p}\bigg(\frac{\binom{2m}{m}}{4^m}\bigg)\bigg|_{m=0}.
%\end{equation}

%\vspace{0.4cm}

We give some examples for specific $p>1$ and $n>0$.

$$\int_{0}^{\pi/2} x\log^2(\sin(x)) \hspace{3pt} dx = \frac{1}{8}\Big(\frac{11\pi^4}{360}+\pi^2\log^22-7\zeta(3)\log2+4\zeta(\overline{3},1)\Big)$$

$$\int_{0}^{\pi/2} x^2\log^2(\sin(x)) \hspace{3pt} dx = \frac{\pi}{24}\Big(\frac{\pi^4}{40}+\pi^2\log^22-9\zeta(3)\log2+12\zeta(\overline{3},1)\Big) $$

\begin{multline*}
\int_{0}^{\pi/2} x^3\log^2(\sin(x)) \hspace{3pt} dx = \frac{1}{64}\Big(\frac{23\pi^6}{420}+\pi^4\log^22+24\pi^2\zeta(\overline{3},1)-48\zeta(\overline{5},1)-24\zeta^2(3)\\-18\pi^2\zeta(3)\log2+93\zeta(5)\log2\Big)
\end{multline*}

\vspace{0.5cm}

Note that if we used $G(n,m,z)$, we would find

\begin{multline}
\frac{\partial^p G}{\partial m^p}\bigg|_{(n,0,\pi)} = - 2^p\Ls_{p+n+1}^{(n)}(\pi) = 2^p\int_{0}^{\pi} x^n\log^{p}\Big(2\sin\Big(\frac{x}{2}\Big)\Big) \hspace{3pt} dx = \frac{\pi^{n+1}}{n+1}\frac{d^p}{dm^p}\binom{2m}{m}\bigg|_{m=0}\\-2n!\sum_{k=1}^{\infty} \frac{\partial^p}{\partial m^p}\binom{2m}{m+k}\bigg|_{m=0}\bigg(\sum_{j=1}^{\lfloor \frac{n+1}{2} \rfloor}\frac{(-1)^j\pi^{n+1-2j}}{(n+1-2j)!k^{2j}}-\delta_{\lfloor \frac{n+1}{2} \rfloor, \frac{n+1}{2}} \frac{(-1)^k(-1)^{\lfloor \frac{n+1}{2}\rfloor}}{k^{n+1}}\bigg).
\end{multline}

\vspace{0.3cm}

%and if $n=0$,

%\begin{equation}
%- 2^p\Ls_{p+1}(\pi) = 2^p\int_{0}^{\pi} \log^{p}\Big(2\sin\Big(\frac{x}{2}\Big)\Big) \hspace{3pt} dx = \frac{\pi}{2}\frac{d^p}{dm^p}\binom{2m}{m}\bigg|_{m=0}.
%\end{equation}

%\vspace{0.5cm}

Again, we compute some integrals below.

$$ \Ls_{4}^{(1)}(\pi) = -\frac{11\pi^4}{720}-2\zeta(\overline{3},1), \hspace{10pt} \Ls_{5}^{(2)} = -\frac{\pi^5}{120}-4\pi\zeta(\overline{3},1)  $$

$$ \Ls_{6}^{(3)}(\pi) = -\frac{23\pi^6}{1680}-6\pi^2\zeta(\overline{3},1)+6\zeta^2(3)+12\zeta(\overline{5},1) $$

$$ \Ls_{7}^{(4)}(\pi) = -\frac{\pi^7}{420}-8\pi^3\zeta(\overline{3},1)+48\pi\zeta(\overline{5},1) $$

\vspace{0.4cm}

\section{Derivatives of binomial coefficients}

\vspace{0.3cm}

These results rely on an efficient calculation of derivatives of central binomial coefficients and shifted central binomial coefficients. In this section, we will provide a proof of a formula for the $p$-th derivative of the binomial coefficients in the above formulae. For simplicity, we will denote these binomial coefficients as if $m$ and $k$ were integers, though the proofs intrinsically use the gamma function (e.g., when taking derivatives). 

\bt{}

Let $ \displaystyle C(m) := \binom{2m}{m+k}$ and $ \displaystyle C_{L}(m) := \frac{1}{4^m}\binom{2m}{m+k}$. Then, we have

\begin{equation}
C^{(p)}(m) = (-1)^p\binom{2m}{m+k}\sum_{i=0}^{p}\binom{p}{i}\Delta_{1}^{p-i}(m)x_{i}(m)
\end{equation}

and

\begin{equation}
C^{(p)}_{L}(m) = (-1)^p\binom{2m}{m+k}\sum_{i=0}^{p}\binom{p}{i}(\Delta_{1}(m)+\log4)^{p-i}x_{i}(m),
\end{equation}

where 

$$x_{0}=1, \hspace{10pt} x_{1}=0, \hspace{10pt} x_{j}(m) = \sum_{l=0}^{j-2} \binom{j-1}{l}\Delta_{j-l}(m)x_{l}(m) $$

\vspace{0.3cm}

with 

$$\Delta_{j}(m) = (-1)^j\big[2^j\psi^{(j-1)}(2m+1)-\psi^{(j-1)}(m+1+k)-\psi^{(j-1)}(m+1-k)\big].$$

\et

\vspace{0.4cm}

\textit{Proof.} We will only provide the proof for $(25)$ as the proof of $(26)$ is identical. The proof is by induction. One can easily see that $C'(m) = -\binom{2m}{m+k}\Delta_{1}(m)$ by expanding out the binomial in terms of the gamma function. Before we move on, we will need a lemma. Further we will omit the argument $m$ throughout the proof.

\bl{For $j \in \mathbb{N}$,}

$$ \frac{d\Delta_{j}}{dm} = -\Delta_{j+1}$$

and

$$ \frac{dx_{j}}{dm} = -x_{j+1}+j\Delta_{2}x_{j-1}.$$

\el

\vspace{0.3cm}

\textit{Proof of Lemma.} The first equation is clear using the definition of $\psi^{(n)}(z)$ and the chain rule. For the second equation, we will do induction. From the recursive definition, $x_{2}=\Delta_{2}$ and so for $j=1$, the second equation is satisfied. Now using the recursive relation for $x_{j}$, the product rule for derivatives, and binomial identities,

$$\frac{dx_{j+1}}{dm} = \frac{d}{dm} \bigg(\sum_{l=0}^{j-1} \binom{j}{l}\Delta_{j+1-l}x_{l}\bigg) = \sum_{l=0}^{j-1}\binom{j}{l}(-\Delta_{j+2-l})x_{l}+\sum_{l=1}^{j-1}\binom{j}{l}\Delta_{j+1-l}(-x_{l+1}+l\Delta_{2}x_{l-1})$$

$$ = -\Delta_{j+2}x_{0}-\sum_{l=1}^{j}\binom{j}{l}\Delta_{j+2-l}x_{l}+\Delta_{2}x_{j}-\sum_{l=0}^{j-1}\binom{j}{l}\Delta_{j+1-l}x_{l+1}+j\Delta_{2}\sum_{l=1}^{j-1}\binom{j-1}{l-1}\Delta_{j-(l-1)}x_{l-1}.$$

\vspace{0.3cm}

Reindexing the second and third sum appropriately, 

$$ \frac{dx_{j+1}}{dm} = -\Delta_{j+2}x_{0}-\sum_{l=1}^{j}\bigg[\binom{j}{l}+\binom{j}{l-1}\bigg]\Delta_{j+2-l}x_{l} + \Delta_{2}x_{j} + j\Delta_{2}x_{j}$$

$$ = -x_{j+2}+(j+1)\Delta_{2}x_{j},$$

which proves this lemma. $\square$

Now going back to the proof of the theorem, by induction we have

\begin{multline*}
\frac{dC^{(p)}}{dm} = (-1)^p\bigg(-\binom{2m}{m+k}\Delta_{1}\bigg)\sum_{i=0}^{p}\binom{p}{i}\Delta_{1}^{p-i}x_{i} + (-1)^p\binom{2m}{m+k}\sum_{i=0}^{p}\binom{p}{i}(p-i)\\ *\Delta_{1}^{p-i-1}(-\Delta_{2})x_{i}+ (-1)^p\binom{2m}{m+k}\sum_{i=1}^{p}\binom{p}{i}\Delta_{1}^{p-i}(-x_{i+1}+i\Delta_{2}x_{i-1})
\end{multline*}

\begin{multline*}
= (-1)^{p+1}\binom{2m}{m+k}\Bigg(\sum_{i=0}^{p+1}\binom{p}{i}\Delta_{1}^{p-i+1}x_{i} + \Delta_{2}\sum_{i=0}^{p}\binom{p}{i}(p-i)\Delta_{1}^{p-i-1}x_{i}\\+\sum_{i=2}^{p+1}\binom{p}{i-1}\Delta_{1}^{p-i+1}x_{i} - \Delta_{2}\sum_{i=1}^{p}\binom{p}{i}i\Delta_{1}^{p-i}x_{i-1}\Bigg)
\end{multline*}

\begin{multline*}
=(-1)^{p+1}\binom{2m}{m+k}\Bigg(\Delta_{1}^{p+1}x_{0}+\sum_{i=1}^{p+1}\bigg[\binom{p}{i}+\binom{p}{i-1}\bigg]\Delta_{1}^{p-i+1}x_{i} \\+ \Delta_{2}\sum_{i=0}^{p-1}\frac{p!}{i!(p-i-1)!}\Delta_{1}^{p-i-1}x_{i} - \Delta_{2}\sum_{i=1}^{p}\frac{p!}{(i-1)!(p-i)!}\Delta_{1}^{p-i}x_{i-1}\Bigg),
\end{multline*}

\vspace{0.3cm}

and now by reindexing, the last two sums cancel. Using the binomial identity as we did in the lemma,

$$ C^{(p+1)}(m) = \frac{dC^{(p)}}{dm} = (-1)^{p+1}\binom{2m}{m+k}\sum_{i=0}^{p+1}\binom{p+1}{i}\Delta_{1}^{p+1-i}x_{i}, $$

\vspace{0.2cm}

which proves the theorem. $\qed$

\vspace{0.2cm}

For simplicity, introduce the following notation:

$$X[s_{n}] := \sum_{j=0}^{n} \binom{n}{j}s_{1}^{n-j}x_{j}$$

where 

$$x_{0}=1, \hspace{10pt} x_{1}=0, \hspace{10pt} x_{j} = \sum_{l=0}^{j-2} \binom{j-1}{l}s_{j-l}x_{l}$$

\vspace{0.3cm}

For completeness, we write out a few sums for an arbitrary sequence $s_{n}$.

$$ X[s_{0}] = 1, \hspace{10pt} X[s_{1}] = s_{1}, \hspace{10pt} X[s_{2}] = s_{1}^2+s_{2} $$

$$ X[s_{3}] = s_{1}^3+3s_{1}s_{2}+s_{3}, \hspace{10pt} X[s_{4}] = s_{1}^4 + 6s_{1}^2s_{2} + 4s_{1}s_{3} + s_{4}+3s_{2}^2 $$

$$ X[s_{5}] = s_{1}^5 + 10s_{1}^3s_{2}+10s_{1}^2s_{3} + 5s_{1}(s_{4}+3s_{2}^2)+s_{5}+10s_{2}s_{3} $$

\vspace{0.4cm}

In fact, these polynomials are known as the complete Bell polynomials, $X[s_{n}] = B_{n}(s_{1}, s_{2}, \dots, s_{n})$ (see \cite{Bouroubi}, \cite{Collins}, \cite{Mihoubi}). Now, evaluating $(25)$ at $m=0$ and letting $\Delta_{j}(0) \equiv \overline{\Delta}_{j} = (-1)^j[2^j\psi^{(j-1)}(1)-\psi^{(j-1)}(1+k)-\psi^{(j-1)}(1-k)]$. Using the reflection formula for the gamma function, we can write $(25)$ as

\begin{equation} 
C^{(p)}(0) = (-1)^p\binom{0}{k}X[\overline{\Delta}_{p}] = (-1)^p\frac{\sin(k\pi)}{k\pi}X[\overline{\Delta}_{p}].
\end{equation}

\vspace{0.3cm}

Using $(9)$ and $(10)$, we can say

\begin{multline*}
\overline{\Delta}_{j} = (-2)^j\psi^{(j-1)}(1) - (-1)^j\Bigg((-1)^{j-1}\pi\frac{d^{j-1}}{dk^{j-1}}\big(\cot(k\pi)\big)\\+(1+(-1)^{j-1})\psi^{(j-1)}(k)+(-1)^{j-1}\frac{(j-1)!}{k^j}\Bigg)
\end{multline*}

$$ = (-2)^j\psi^{(j-1)}(1)+(1-(-1)^j)\psi^{(j-1)}(k)+\frac{(j-1)!}{k^j}+\frac{d^{j-1}}{dk^{j-1}}\big(\pi\cot(k\pi)\big) \equiv \xi_{j}+ \nu_{j}, $$

\vspace{0.3cm}

where

$$ \nu_{j} := \frac{d^{j-1}}{dk^{j-1}}\big(\pi\cot(k\pi)\big), \hspace{15pt} \xi_{j} := \overline{\Delta}_{j}-\nu_{j}. $$

\vspace{0.3cm}

Thus we have

$$C^{(p)}(0) = (-1)^p\frac{\sin(k\pi)}{k\pi}X[\xi_{p}+\nu_{p}].$$

\vspace{0.3cm}

Using the binomial theorem and some algebra, one can see

\begin{equation}
C^{(p)}(0) = (-1)^p\frac{\sin(k\pi)}{k\pi}\sum_{i=0}^{p}\binom{p}{i} X[\xi_{p-i}]X[\nu_{i}],
\end{equation}

\vspace{0.3cm}

which is also a well-known binomial identity of the complete Bell polynomials. We can simplify $(28)$ more with the help of two lemmas.

\vspace{0.4cm}

\bl For $n \in \mathbb{N}_{0}$,

$$X[\nu_{n}] =  \left\{
\begin{array}{ll}
      (-1)^j\pi^{2j} & n = 2j, \hspace{5pt} j \in \mathbb{N}_{0}, \\ \\
      (-1)^j\pi^{2j+1}\cot(k\pi) & n=2j+1, \hspace{5pt} j \in \mathbb{N}_{0}.\\
\end{array} 
\right.
$$

\el

\textit{Proof.} Note that 

$$ \nu_{n}=\frac{d^{n-1}}{dk^{n-1}}\big(\pi\cot(k\pi)\big) = \frac{d^n}{dk^n}\big(\log(\sin(k\pi))\big).$$

\vspace{0.3cm}

So, letting $y=\log(\sin(k\pi))$ and using the formula provided in \cite{Collins}, 

$$ B_{n}(\nu_{1}, \nu_{2}, \dots, \nu_{n}) = X[\nu_{n}] = e^{-y}\frac{d^ne^y}{dk^n} = \csc(k\pi)\frac{d^n}{dk^n}\big(\sin(k\pi)\big). $$

\vspace{0.3cm}

When $n=2j$, we see $\displaystyle \frac{d^{2j}}{dk^{2j}}\big(\sin(k\pi)\big) = (-1)^j\pi^{2j}\sin(k\pi)$ and when $n=2j+1$, we have $\displaystyle \frac{d^{2j+1}}{dk^{2j+1}}\big(\sin(k\pi)\big) = (-1)^j\pi^{2j+1}\cos(k\pi)$. So the proof is complete. $\qed$

\vspace{0.3cm}

Since our results involve a sum over natural numbers $k$, we can simplify $(28)$ to

$$ C^{(p)}(0) = (-1)^p\frac{\sin(k\pi)}{k\pi}\sum_{j=0}^{\lfloor \frac{p-1}{2} \rfloor}
\binom{p}{2j+1} X[\xi_{p-1-2j}](-1)^j\pi^{2j+1}\cot(k\pi) $$

$$ = \frac{(-1)^{p+k}}{k}\sum_{j=0}^{\lfloor \frac{p-1}{2} \rfloor}\binom{p}{2j+1}X[\xi_{p-1-2j}](-\pi^2)^j. $$

\vspace{0.3cm}

\bl Let $\rho_{1} = \pi^2/3$, and 

\vspace{0.2cm}

$$ \rho_{n} = (-1)^{n+1}\Bigg(\frac{\pi^{2n}}{2n+1}+\sum_{i=1}^{n-1}\binom{2n-1}{2i-1}\frac{(-1)^i\pi^{2n-2i}}{2n-2i+1}\rho_{i}\Bigg).$$

\vspace{0.2cm}

Then, $\rho_{n} = 2\psi^{(2n-1)}(1)$.

\el

\vspace{0.4cm}

\textit{Proof.} The proof is by induction. It is clearly true for $n=1$, so using our induction hypothesis, along with $(5)$ and $(11)$,

$$ \rho_{n} = (-1)^{n+1}\Bigg(\frac{\pi^{2n}}{2n+1}+\sum_{i=1}^{n-1}\binom{2n-1}{2i-1}\frac{(-1)^i\pi^{2n-2i}}{2n-2i+1}\rho_{i}\Bigg)$$

$$ = (-1)^{n+1}\Bigg(\frac{\pi^{2n}}{2n+1}+\sum_{i=1}^{n-1}\binom{2n-1}{2i-1}\frac{(-1)^i\pi^{2n-2i}}{2n-2i+1}(2i-1)!\frac{(-1)^{i+1}B_{2i}(2\pi)^{2i}}{(2i)!}\Bigg)$$

$$ = (-1)^{n+1}\Bigg(\frac{\pi^{2n}}{2n+1}-\sum_{i=1}^{n-1}\frac{(2n-1)!}{(2n-2i)!}\frac{\pi^{2n}}{2n-2i+1}\frac{B_{2i}2^{2i}}{(2i)!}\Bigg)$$

$$ = (-1)^{n+1}\Bigg(\frac{\pi^{2n}}{2n+1}-\frac{\pi^{2n}}{2n(2n+1)}\sum_{i=1}^{n-1}\binom{2n+1}{2i}B_{2i}2^{2i}\Bigg).$$

\vspace{0.3cm}

Using the fact that $B_{0}=1$, $B_{1}=-1/2$, and $B_{2k+1}=0$ for $k \in \mathbb{N}$, we can rewrite the sum and obtain

$$\rho_{n} = \frac{(-1)^{n+1}\pi^{2n}}{2n}\Bigg(\frac{2n}{2n+1}-\frac{1}{2n+1}\bigg(-1+(2n+1)-(2n+1)B_{2n}2^{2n}+\sum_{i=0}^{2n+1} \binom{2n+1}{i}B_{i}2^i\bigg)\Bigg)$$

$$ = \frac{(-1)^{n+1}\pi^{2n}}{2n}\Bigg(B_{2n}2^{2n}-\frac{2^{2n+1}}{2n+1}\sum_{i=0}^{2n+1} \binom{2n+1}{i}B_{i}\Big(\frac{1}{2}\Big)^{2n+1-i}\Bigg),$$

\vspace{0.3cm}

and using the identities $\displaystyle \sum_{i=0}^{k} \binom{k}{i}B_{i}z^{i} = B_{k}(z)$ where $B_{k}(z)$ is the $k$-th Bernoulli polynomial, and $\displaystyle B_{k}\Big(\frac{1}{2}\Big)=\bigg(\frac{1}{2^{k-1}}-1\bigg)B_{k}$, this sum completely vanishes. So this simplifies to

$$\rho_{n} = \frac{(-1)^{n+1}\pi^{2n}}{2n}\bigg(B_{2n}2^{2n}-0\bigg) = 2\psi^{(2n-1)}(1), $$

using $(5)$ and $(11)$ again. $\qed$

%This formula does not appear to be in any literature either and has been confirmed up to $n=20$. This formula also immediately implies a recursive formula for $\zeta(2n)$, the even arguments of the Riemann zeta function. Using conjecture 1 and conjecture 2, we are able to provide the main equation for computing derivatives of shifted binomial coefficients. 

\vspace{0.3cm}

Now we can state the second main theorem of this section.

\bt For $k$ in $\mathbb{N}$,

\begin{equation}
C^{(p)}(0) = \frac{(-1)^{p+k}}{k}pX[\overline{\xi}_{p-1}] = \frac{(-1)^{p+k}}{k}pB_{p-1}(\overline{\xi}_{1},\overline{\xi}_{2},\dots,\overline{\xi}_{p-1})
\end{equation}

\vspace{0.3cm}

where

$$\overline{\xi}_{j} = \xi_{j}-(1+(-1)^j)\psi^{(j-1)}(1), $$

that is,

$$ \overline{\xi}_{j} = (-2)^j\psi^{(j-1)}(1) + 2\delta_{\lfloor \frac{j+1}{2} \rfloor, \frac{j+1}{2}}\psi^{(j-1)}(k)-2\delta_{\lfloor \frac{j}{2} \rfloor, \frac{j}{2}}\psi^{(j-1)}(1) + \frac{(j-1)!}{k^j}. $$

\et

\vspace{0.3cm}

\textit{Proof.} To prove this, notice

$$ \frac{(-1)^{p+k}}{k}pX[\overline{\xi}_{p-1}] = \frac{(-1)^{p+k}}{k}pX[\xi_{p-1}-(1+(-1)^{p-1})\psi^{(p-2)}(1)]$$

$$ = \frac{(-1)^{p+k}}{k}pX[\xi_{p-1}+\alpha_{p-1}] = \frac{(-1)^{p+k}}{k}p\sum_{i=0}^{p-1}\binom{p-1}{i}X[\xi_{p-1-i}]X[\alpha_{i}], $$

\vspace{0.3cm}

where $\alpha_{i} = -(1+(-1)^i)\psi^{(i-1)}(1)$, that is, $\alpha_{2i} = -\rho_{i}$ and $\alpha_{2i+1}=0$. In particular, $\alpha_{1}=0$ so from the definition of $X$, $X[\alpha_{i}] = x_{i}$, where 

$$x_{0} = 1, \hspace{10pt} x_{1}=0, \hspace{10pt} x_{i} = \sum_{l=0}^{i-2} \binom{i-1}{l}\alpha_{i-l}x_{l}.$$
\vspace{0.1cm}

\textit{Claim.} For $i \in \mathbb{N}_{0}$,

$$x_{2i+1} = 0, \hspace{15pt} x_{2i} = \frac{(-1)^i\pi^{2i}}{2i+1}. $$

\vspace{0.2cm}

\textit{Proof of Claim.} $x_{1}=0$ by definition, so assume $x_{k}=0$ for odd $k \leq 2i-1$. We can see that $x_{2i+1}$ will be a sum of $\alpha_{2i+1-l}$ and $x_{l}$ multiplied together for $l \leq 2i-1$. If $l$ is even, $2i+1-l$ is odd and so $\alpha_{2i+1-l}=0$. If $l$ is odd, by the induction assumption, $x_{l}=0$. So all terms of the sum will be 0 and thus $x_{2i+1}=0$ for all $i \in \mathbb{N}_{0}$. For the even indices, we will also use induction. For $i=0$, this is clearly satisfied. Now, assume the formula for indices less than $2i$. Using $x_{2i+1}=0$ along with the definition of $\alpha_{k}$,

$$ x_{2i} = \sum_{l=0}^{2i-2} \binom{2i-1}{l}\alpha_{2i-l}x_{l} = \alpha_{2i} + \sum_{l=2}^{2i-2}\binom{2i-1}{l}\alpha_{2i-l}x_{l}$$

$$ = \alpha_{2i}+\sum_{l=1}^{i-1}\binom{2i-1}{2l}\alpha_{2i-2l}x_{2l} = -\rho_{i}-\sum_{l=1}^{i-1}\binom{2i-1}{2(i-l)-1}\rho_{i-l}x_{2l}. $$

\vspace{0.3cm}

Using a change of index on the sum, our induction hypothesis, and the previous lemma about $\rho_{n}$,

$$ x_{2i} = -\rho_{i}-\sum_{l=1}^{i-1} \binom{2i-1}{2l-1}\rho_{l}x_{2i-2l} = -\rho_{i}-\sum_{l=1}^{i-1}\binom{2i-1}{2l-1}\frac{(-1)^{i-l}\pi^{2i-2l}}{2i-2l+1}\rho_{l} $$

$$ = -\rho_{i} - (-1)^i\bigg((-1)^{i+1}\rho_{i}-\frac{\pi^{2i}}{2i+1}\bigg) = \frac{(-1)^i\pi^{2i}}{2i+1},$$

\vspace{0.3cm}

and so the claim is proven. $\square$ 

Now, we are able to write

$$ \frac{(-1)^{p+k}}{k}pX[\overline{\xi}_{p-1}] = \frac{(-1)^{p+k}}{k}p\sum_{i=0}^{p-1}\binom{p-1}{i}X[\xi_{p-1-i}]x_{i}$$

$$ = \frac{(-1)^{p+k}}{k}p\sum_{j=0}^{\lfloor \frac{p-1}{2}\rfloor} \binom{p-1}{2j}X[\xi_{p-1-2j}]\frac{(-1)^j\pi^{2j}}{2j+1} $$

$$= \frac{(-1)^{p+k}}{k}\sum_{j=0}^{\lfloor \frac{p-1}{2} \rfloor}\binom{p}{2j+1}X[\xi_{p-1-2j}](-\pi^2)^j = C^{(p)}(0), $$

\vspace{0.3cm}

which proves the theorem. $\qed$

\vspace{0.3cm}

Note that the two lemmas also imply a very similar formula for $(26)$, the only difference is $\overline{\xi}_{L,j} = \overline{\xi}_{j}+\delta_{j1}\log4 $. For completeness we write some of these out, in their simplified form where $H_{k}$ is the $k$-th harmonic number.

\vspace{0.3cm}

$$ \frac{d}{dm}\binom{2m}{m+k}\bigg|_{m=0} = \frac{(-1)^{k+1}}{k}$$

$$ \frac{d^2}{dm^2}\binom{2m}{m+k}\bigg|_{m=0} = \frac{2(-1)^k}{k}\Big(2H_{k}-\frac{1}{k}\Big) $$

$$ \frac{d^3}{dm^3}\binom{2m}{m+k}\bigg|_{m=0} = \frac{3(-1)^{k+1}}{k}\Bigg(\Big(2H_{k}-\frac{1}{k}\Big)^2+\frac{1}{k^2}+2\psi^{(1)}(1)\Bigg)$$

\begin{multline*}
\frac{d^4}{dm^4}\binom{2m}{m+k}\bigg|_{m=0} = \frac{4(-1)^k}{k}\Bigg(\Big(2H_{k}-\frac{1}{k}\Big)^3+3\Big(2H_{k}-\frac{1}{k}\Big)\\ * \Big(2\psi^{(1)}(1)+\frac{1}{k^2}\Big)+2\psi^{(2)}(k)+\frac{2}{k^3}-8\psi^{(2)}(1)\Bigg)
\end{multline*}

\vspace{0.4cm}

Lastly, for the central binomial coefficients, i.e., when $k=0$, we can still use equation $(25)$. Let

$$\eta_{j}(m) := \Delta_{j}(m)\Big|_{k=0} = (-1)^j(2^j\psi^{(j-1)}(2m+1)-2\psi^{(j-1)}(m+1)).$$

\vspace{0.5cm}

Then we have

$$\frac{d^p}{dm^p}\binom{2m}{m} = (-1)^p\binom{2m}{m}X[\eta_{p}].$$

\vspace{0.5cm}

Letting $\overline{\eta}_{j} \equiv \eta_{j}(0) = (-1)^j(2^j-2)\psi^{(j-1)}(1)$, note that $\overline{\eta}_{1}=\eta_{1}(0)=0$ so all terms vanish except the $j=p$ term in the definition of $X$. So for $m=0$,

\begin{equation}
\frac{d^p}{dm^p}\binom{2m}{m}\Bigg|_{m=0} = (-1)^px_{p} = (-1)^pB_{p}(0,\overline{\eta}_{2},\overline{\eta}_{3},\dots,\overline{\eta}_{p})
\end{equation}

and

\begin{equation}
\frac{d^p}{dm^p}\Bigg(\frac{1}{4^m}\binom{2m}{m}\Bigg)\Bigg|_{m=0} = (-1)^p\sum_{i=0}^{p}\binom{p}{i}(\log4)^{p-i}x_{i},
\end{equation}

\vspace{0.3cm}

where

$$ x_{0} = 1, \hspace{10pt} x_{1} = 0, \hspace{10pt} x_{n} = \sum_{l=0}^{n-2}\binom{n-1}{l}\overline{\eta}_{n-l}x_{l}. $$

\vspace{0.5cm}

Using $(30)$ and our previous results, we have

\begin{equation}
\Ls_{p+n+1}^{(n)}(z) = -\int_{0}^{z} x^n\log^p\Big(2\sin\Big(\frac{x}{2}\Big)\Big) \hspace{3pt} dx =\frac{(-1)^{p+1}z^{n+1}}{2^p(n+1)}B_{p}(0,\overline{\eta}_{2},\overline{\eta}_{3},\dots,\overline{\eta}_{p})
\end{equation}

\vspace{0.3cm}

for $z=2\pi, n\in \{0,1\}$ or $z=\pi,n=0$. For $n=0$ and $z\in \{\pi,2\pi\}$, this formula is well-known (see \cite{Borwein}, \cite{Choi}). For $n=0$ and arbitrary $z$,

\begin{multline}
\Ls_{p+1}(z)=-\int_{0}^{z} \log^p\Big(2\sin\Big(\frac{x}{2}\Big)\Big) \hspace{3pt} dx = \frac{(-1)^{p+1}}{2^p}\bigg(zB_{p}(0,\overline{\eta}_{2},\overline{\eta}_{3},\dots,\overline{\eta}_{p})\\+2p\sum_{k=1}^{\infty}\frac{\sin(kz)}{k^2}B_{p-1}(\overline{\xi}_{1},\overline{\xi}_{2},\dots,\overline{\xi}_{p-1})\bigg),
\end{multline}

\vspace{0.3cm}

where 

$$\overline{\eta}_{j} = (-1)^j(2^j-2)\psi^{(j-1)}(1),$$

and

$$ \overline{\xi}_{j} = (-2)^j\psi^{(j-1)}(1) + 2\delta_{\lfloor \frac{j+1}{2} \rfloor, \frac{j+1}{2}}\psi^{(j-1)}(k)-2\delta_{\lfloor \frac{j}{2} \rfloor, \frac{j}{2}}\psi^{(j-1)}(1) + \frac{(j-1)!}{k^j} $$

$$= \overline{\eta}_{j}+2\delta_{\lfloor \frac{j+1}{2} \rfloor, \frac{j+1}{2}}(\psi^{(j-1)}(k)-\psi^{(j-1)}(1)) + \frac{(j-1)!}{k^j}. $$

\bigskip

\begin{flushright}
\begin{minipage}{148mm}\sc\footnotesize
University of Pittsburgh, Department of Mathematics, 301 Thackeray Hall, Pittsburgh, PA 15260, USA\\
{\it E--mail address}: {\tt djo15@pitt.edu} \vspace*{3mm}
\end{minipage}
\end{flushright}

\end{document}